\documentclass[a4wide,10pt]
{article}
\usepackage{amssymb, amsthm,mathtools,
verbatim, amsmath}
\usepackage[arrow, matrix, curve]{xy}
\usepackage{hyperref}
\usepackage{mathrsfs}

\RequirePackage[dvipsnames,usenames]{color}  

\usepackage{mabliautoref}

\input{GBplainNew2.sty}

\long\def\forget#1{}

\advance\textheight by 1in
\advance\topmargin by -.4in

\maketheorem{Thm}{Theorem}{theorem}
\maketheorem{Prop}{Proposition}{theorem}
\maketheorem{Prop-Def}{Proposition-Definition}{theorem}
\maketheorem{Cor}{Corollary}{theorem}
\maketheorem{Lem}{Lemma}{theorem}
\maketheorem{SubLem}{Sublemma}{theorem}
\maketheorem{Conj}{Conjecture}{theorem}

\theoremstyle{definition}    

\maketheorem{Def}{Definition}{theorem}
\maketheorem{DefRem}{Definition-Remark}{theorem}
\maketheorem{Claim}{Claim}{theorem}
\maketheorem{Alg}{Algorithm}{theorem}

\theoremstyle{remark}    

\maketheorem{Rem}{Remark}{theorem}
\maketheorem{Ex}{Example}{theorem}
\maketheorem{Ques}{Question}{theorem}
\maketheorem{QuesRem}{Question-Remark}{theorem}
\maketheorem{Ass}{Assumption}{theorem}
\maketheorem{Con}{Convention}{theorem}

\forget{\newif\ifnormalesBeweisEnde

  {\vskip 0.3ex plus 0.5ex minus 0ex \pagebreak[1]
   \global\normalesBeweisEndetrue
   \trivlist
   \item[\hskip\labelsep {\textsc{Proof} \rm #1}:]}%
  {\ifnormalesBeweisEnde \EndOfProof \fi
   \endtrivlist
   \vskip 1ex plus 1ex minus 0ex \pagebreak[2]}
  {\vskip 0.3ex plus 0.5ex minus 0ex \pagebreak[1]
   \global\normalesBeweisEndetrue
   \trivlist
   \item[\hskip\labelsep \textsc{Proof}:]}%
  {\ifnormalesBeweisEnde \EndOfProof \fi
   \endtrivlist
   \vskip 1ex plus 1ex minus 0ex \pagebreak[2]}
\def\EndOfProof{\hskip .5em \vrule width 1.0ex height 1.0ex depth 0.3ex}
}

\usepackage{latexsym}
\usepackage{amsmath}
\usepackage{amstext}
\usepackage{amssymb}
\usepackage{amsbsy}

\def\theenumi{(\alph{enumi})}

\def\theenumii{(\roman{enumii})}
\def\p@enumii{\theenumi}

\newdimen\auxdimen
\newdimen\auxdim
\newdimen\auxdimone
\newdimen\auxdimtwo
\newdimen\auxdimthree

\long\def\forget#1{}

\def\?{\ ???\ \immediate\write16{}%
\immediate\write16{Warning: There was still a question mark . . . }%
\immediate\write16{}}

\newcommand{\X}{x}
\newcommand{\Z}{z}
\newcommand{\Y}{y}
\newcommand{\spec}{\mathrm{spec}}
\renewcommand{\et}{{\acute{\mathrm{e}}\mathrm{t}}}

\DeclareMathOperator{\CNL}{{CNL}_{\breve \CO}}
\DeclareMathOperator{\Tate}{Ta}
\DeclareMathOperator{\Spf}{Spf}
\newcommand{\tempnewpage}{}%

\renewcommand{\to}[1][]{\xrightarrow{\ #1\ }}

\advance\textwidth by 1in
\advance\oddsidemargin by -.5in

\title{The generic monodromy of Drinfeld modular varieties\\ in special characteristic}
\author{Gebhard B\"ockle and Florian Breuer}

\begin{document}
\parindent0em
\parskip.3em

\maketitle

\begin{abstract}
By combining theorems of Drinfeld and Strauch, we show that the monodromy representation on the special fibre of a Drinfeld modular variety, with level not divisible by the characteristic, is surjective. We illustrate this result in the special case of Drinfeld $\BF_q[t]$-modules in level $t$, and apply this to show that the Kronecker factors of a Drinfeld modular polynomial in rank $r$ are irreducible.
\end{abstract}

\begin{center}
{\em
Dedicated to Gerhard Frey on the occasion of his $75^{\mathrm{th}}$ birthday.
}
\end{center}


\parindent0em
\parskip.3em

\tempnewpage

\section{Statement of the main result}
\label{Sect:One}

Throughout this paper we fix a global function field $F$ of characteristic $p>0$ with exact field of constants the finite field $\BF_q$ of cardinality $q$. We fix a place $\infty$ of $F$, and let $A$ denote the ring of elements of $F$ which are regular away from $\infty$. This is a Dedekind domain with finite class group $\Cl(A)$ and unit group $A^\times=\BF_q^\times$. 

Let $I\subset A$ denote a proper non-zero ideal and $n_I$ the order of $I$ in the ideal class group $\Cl(A)$ of $A$. Let $g_I\in A$ be a generator of $I^{n_I}$; it is unique up to multiplication by $\BF_q^\times$. Hence $A[1/I]:=A[1/g_I]$ is independent of the choice of~$g_I$. We also write $\Cl_I(A)$ for the $I$-class group of $A$, i.e., the group of fractional ideals of $A$ of support prime to $I$ modulo its subgroup of principal fractional ideals that possess a generator which is congruent to $1$ modulo~$I$. One has a short exact sequence $0\to (A/I)^\times/\BF_q^\times\to \Cl_I(A)\to \Cl(A)\to 0$.

Consider the functor $\CM_I^r$ from $A[1/I]$-schemes $S$ to $\Sets$, which to any such $S$ assigns the set of isomorphism classes of tuples $(\CL,\phi,\alpha)$, where $\CL$ is a line bundle on $S$, where $\phi$ (together with $\CL$) is a Drinfeld $A$-module $\phi\colon A\to\End_{\BF_q\textrm{-gp sch.}/S}(\CL)$ 
of rank $r$, and where $\alpha$ denotes a level $I$-structure on $(\CL,\phi)$, subject to the condition that the characteristic $\partial \phi\colon A\to \End_S(\Lie \CL)=\CO_S$ coincides with the structure morphism $S\to\Spec A[1/I]$ composed with the open immersion $\Spec A[1/I]\to\Spec A$. By \cite[Prop.~5.3 and Cor.~to 5.4]{DrinfeldEM}, the functor $\CM_I^r$ is representable by a smooth finite type morphism $\FM_I^r\to\Spec A[1/I]$ of relative dimension $r-1$. The universal Drinfeld module on $\FM_I^r$ we denote by
\[\phi_I^r\colon A\to \End_{\BF_q\textrm{-gp sch.}/\FM_I^r}(\CL_I^r).\]

Let now $\Fp\subset A$ denote a maximal ideal that is prime to $I$. We write $\kappa_\Fp$ for its residue field, 
 $A_\Fp$ for the completion of $A$ at $\Fp$, and let $\overline\kappa_\Fp$ be an algebraic closure of $\kappa_\Fp$.

\begin{Def}
We call $\FM_{I,\Fp}^r:=\FM_I^r\times_{\Spec A[1/I]}\Spec \kappa_\Fp$ the \emph{special fiber of $\FM_I^r$ at $\Fp$.}
\end{Def}

Let  $\FM_{I,\overline\Fp}^r$ be the base change $\FM_{I,\Fp}^r\times_{\kappa_\Fp}\overline \kappa_\Fp$ and let $\phi_{I,\overline\Fp}^r$ be the corresponding universal Drinfeld module
. The scheme $\FM_{I,\overline\Fp}^r$  is regular.
Its connected components can be naturally labelled by $\Cl_I(A)$: By \cite[Cor. 4.6]{Papikian} the connected components of $\FM_{I,\overline\Fp}^r$ are in bijection with those of $\FM_{I,\overline\Fp}^1 \cong \Spec \big( R_I\otimes_{A[1/I]}\overline\kappa_{\Fp}\big)$, where by \cite[Thm. 1, \S7]{DrinfeldEM} $R_I$ is the integral closure of $A$ in the class field of $F$ associated to $\Cl_I(A)$. Now class field theory gives the desired labeling.

 For each ideal class $\Fc$ in $\Cl_I(A)$, denote by $\eta_\Fc=\Spec k_{\eta_\Fc}$ the generic point of the corresponding component, and let $\overline\eta_{\Fc}=\Spec \overline \kappa_{\eta_\Fc}$ be a geometric point above $\kappa_{\eta_\Fc}$. Observe that $\kappa_{\eta_\Fc}$ contains~$\overline\kappa_\Fp$.

Let $\phi_{\eta_\Fc}^r$ denote the pullback of $\phi_{I,\overline\Fp}^r$  to $\eta_\Fc$. By \cite[Prop.~5.5]{DrinfeldEM}, it is a Drinfeld $A$-module of characteristic $\Fp$ and height $1$, i.e., $\phi_{\eta_\Fc}^r$ is ordinary. This means that for any $n\ge1$ the group of $\Fp^n$-torsion points $\phi_{\eta_\Fc}^r[\Fp^n](\overline \kappa_{\eta_\Fc})$ is a free $A/\Fp^n$-module of rank $r-1$. 

Let $g_\Fp$ be a generator of the principal ideal $\Fp^{n_\Fp}$, so that $\phi_{\eta_\Fc}^r[\Fp^{nn_\Fp}](\overline \kappa_{\eta_\Fc})$ is the set of roots of $\phi^r_{\eta_\Fc,g_\Fp^n}(X)$. We define the $\Fp$-adic Tate module of $\phi_{\eta_\Fc}$ as 
\[\Tate_\Fp\phi^r_{\eta_\Fc}=\invlim_n \phi_{\eta_\Fc}^r[\Fp^{nn_\Fp}](\overline \kappa_{\eta_\Fc}), \]
where multiplication by $g_\Fp$ defines the transition map in the inverse system. The limit is independent of the choice of $g_\Fp$. By ordinariness of $\phi_{\eta_\Fc}$ it is free of rank $r-1$ over~$A_\Fp$. 

Observe that $\phi^r_{\eta_\Fc,g_\Fp^n}(X)=h_n\circ (X\mapsto X^{q^{n\deg \Fp}})$ for some unique $\BF_q$-linear Polynomial $h_n\in\kappa_{\eta_\Fc}[X]$ with non-vanishing linear term. The \'etale quotient $\phi_{\eta_\Fc}^r[\Fp^{nn_\Fp}]^\et$ of the finite flat $A$-module scheme $\phi_{\eta_\Fc}^r[\Fp^{nn_\Fp}]$ is $\Spec \kappa_{\eta_\Fc}[X]/(h_n(X))$. The group schemes $\phi_{\eta_\Fc}^r[\Fp^{nn_\Fp}]^\et$ also form an inverse system, and one has $\phi_{\eta_\Fc}^r[\Fp^{nn_\Fp}](\overline \kappa_{\eta_\Fc})\cong\phi_{\eta_\Fc}^r[\Fp^{nn_\Fp}]^\et(\kappa^\sep_{\eta_\Fc}) $ as finite $A$-modules. Because the polynomials $h_n$ are defined over $\kappa_{\eta_\Fc}$, the absolute Galois group $G_{\kappa_{\eta_\Fc}}=\Gal(\kappa^\sep_{\eta_\Fc}/\kappa_{\eta_\Fc})$ acts on $\Tate_\Fp\phi^r_{\eta_\Fc}$ and by the very construction of $\Tate_\Fp\phi^r_{\eta_\Fc}$, the group $G_{\kappa_{\eta_\Fc}}$ acts $A_\Fp$-linearly. This yields a continuous group homomorphism
\[ \rho_{\Fp,\eta_\Fc}\colon G_{\kappa_{\eta_\Fc}}\longto \Aut_{A_\Fp}(\Tate_\Fp\phi^r_{\eta_\Fc})\cong\GL_{r-1}(A_\Fp).\]

Our main result is the following:
\begin{Thm}\label{Thm:Main}
The map $\rho_{\Fp,\eta_\Fc}$ is surjective.
\end{Thm}

The proof is a simple consequence of the results \cite{DrinfeldEM,Strauch} by Drinfeld and Strauch, which seems not to have been recorded in the literature. In fact, combining the work of Drinfeld and Strauch, it even follows that the image under $\rho_{\Fp,\eta_\Fc}$ of a decomposition group of $G_{\kappa_{\eta_\Fc}}$ at a supersingular point of $\FM_{I,\overline\Fp}^r$ in the component of $\eta_\Fc$ is already surjective; cf.~\autoref{Rem:DecompGp}.

\medskip

One might wonder about refinements of \autoref{Thm:Main}. For any point $x$ of $\FM_{I,\overline\Fp}^r$ denote by $\rho_{\Fp,x}\colon G_x\to \Aut_{\Fp}(\Tate_\Fp\phi^r_x)$ the action of the absolute Galois group $G_x=\Gal(\kappa_x^\sep/\kappa_x)$ of the residue field at $x$ on the corresponding Tate module $\Tate_\Fp\phi^r_x$ of rank $r_x$ with $0\le r_x\le r-1$. 
\begin{Ques}
Suppose that $\End_{\overline \kappa_x}(\phi^r_x)=A$. Is $\rho_{\Fp,x}(G_x)$ open in $\Aut_{\Fp}(\Tate_\Fp\phi^r_x)\cong \GL_{r_x}(A_\Fp)$? 
\end{Ques}
This appears to be a natural analog of the results \cite{DevicPink} of Devic and Pink on adelic openness for Drinfeld modules in special characteristic. They consider the Galois action of a Drinfeld $A$-module $\phi$ of rank $r$ and characteristic $\Fp\neq0$ over a finitely generated field. If $\End_{\overline \kappa_x}(\phi_x)=A$, their results imply that the associated adelic Galois representation of $G_x$ away from $\Fp$ and $\infty$ has open image in $\SL_r(\prod_{v\neq\Fp}A_v)$. They also give a complete answer with no condition on $\End_{\overline \kappa_x}(\phi_x)$. This leads to.
\begin{Ques}
Describe for any point $x$ of $\FM_{I,\overline\Fp}^d$ the Zariski closure $\CG_x$ of $\rho_{\Fp,x}(G_x)$ in $\Aut_{\Fp}(\Tate_\Fp\phi^r_x)\cong \GL_{r_x}(A_\Fp)$. Is  $\rho_{\Fp,x}(G_x)$ an open subgroup if $\CG_x(A_\Fp)$?
\end{Ques}

We end this introduction with a quick survey of the content of the individual sections. Section 2 recalls the relevant work of Drinfeld on formal $\CO$-modules and $\CO$-divisible groups from~\cite{DrinfeldEM}. Section 3 recalls the main theorem of Strauch, so that in Section 4 we can combine the two and deduce the proof of \autoref{Thm:Main}. Section 5 illustrates the main result in the special case of $A=\BF_q[t]$ and level $t$, where the moduli scheme can be described explicitly. In section 6 we shall answer in \autoref{Prop:OnBR16} a question raised in \cite{BR16} related to the reduction of modular polynomials of $\FM_{\Fp}^r$ in the case $A=\BF_q[t]$. We shall prove that certain special polynomials which are the natural building blocks of the mod $\Fp$ reduction of modular polynomials are irreducible as asked in~\cite[Question~4.5]{BR16}.

\tempnewpage

\section{Formal {$\CO$}-modules, {$\CO$}-divisible groups and deformations of Drinfeld modules}
\label{Sec:Two}

Let $K$ be a non-archimedean local field with ring of integers $\CO$ and finite residue field $k$. The normalized valuation on $K$ is $v_K$, its uniformizer $\varpi_K$ and the cardinality of $k$ will be $q_K$. Let $\breve K$ be the completion of the maximal unramified extension of $K$ and write $\breve\CO$ for its ring of integers. The residue field $\breve k$ of $\breve\CO$ is an algebraic closure of $k$. Denote by $\CNL$ the category of complete noetherian local $\breve\CO$-algebras $C$ with residue field $\breve k$, and with morphisms being $\breve\CO$-algebra homomorphisms $f\colon C\to C'$ such that $f(\Fm_C)\subset\Fm_{C'}$, where for $C\in\CNL$ we denote by $\Fm_C$ its maximal ideal.

Let $B$ be a ring. 
The power series ring over $B$ in indeterminates $\X_1,.\,.\,,\X_n$ will be $B[[\X_1,.\,.\,,\X_n]]$.
\begin{Def}[{\cite[\S~1]{DrinfeldEM}}]
A formal group\footnote{More correctly we should add the attributes one-dimensional and commutative; but we shall not deal with any other kind of formal group; and so for the sake of brevity we suppress them.} over $B$ is a series $\Phi$ in $B[[\X,\Y]]$ such that $\Phi(\X,\Y)=\Phi(\Y,\X)$, $\Phi(\X,0)=\X$ and \hbox{$\Phi(\Phi(\X,\Y),\Z)=\Phi(\X,\Phi(\Y,\Z))$.}

A {\em homomorphism} from a formal group $\Phi$ to a formal group $\Psi$ over $B$ is a series $\beta\in \X B[[\X]]$ such that $\Psi(\beta(\X),\beta(\Y))=\beta(\Phi(\X,\Y))$. Composition of homomorphisms is composition of formal power series; it is well-defined because the series have zero constant term.

The endomorphism  ring of a formal group $\Phi$ is denoted $\End(\Phi)$. It comes with a natural homomorphism $D\colon \End(\Phi)\to B$, given by differentiation at zero $\beta\mapsto D\beta= \big(\frac{\mathrm{d}}{\mathrm{d}\X}\beta\big)(0)$.

Suppose that $B$ is an $\CO$-algebra via a map $\alpha$. A {\em formal $\CO$-module over $B$} is a pair $X=(\Phi,[\cdot]_X)$ where $\Phi$ is a formal group over $B$ and $[\cdot]_X$ is a homomorphism $\CO\to\End(\Phi)$ such that $D\circ[\cdot]_X=\alpha$. Morphisms of formal $\CO$-modules are defined in the obvious way.

If $\alpha(\varpi_K)=0$, then $[\varpi_K]_X\in xB[[\X]]$ 
can be written as the composition $\gamma\circ (\X\mapsto \X^{q_K^h})$ for some unique $\gamma\in \X B[[\X]]$ with linear term $D\gamma\neq0$ and a unique $h\ge 1$. One calls $h$ the {\em height} of the
formal $\CO$-module $X=(\Phi,[\cdot]_X)$.

%
%
%
%
%
 %
\end{Def}

\begin{Ex}[{\cite[Rem.~after~Prop.~2.2]{DrinfeldEM}, \cite[\S~4]{Rosen}}]\label{Ex:DM->FormalGp}
Let $s\colon A\to \CO$ be a ring homomorphism with $s(\Fp)\CO=\varpi_K\CO$, 
and let $C$ be in $\CNL$. This induces an $A$-algebra  structure on $C$, which we denote by $\gamma$. Let $\phi\colon A\to C\{\tau\},a\mapsto \phi_a$ be a Drinfeld $A$-module in standard form of rank $r$ and characteristic $\gamma$; cf.~\cite[Rem.~after Prop.~5.2]{DrinfeldEM}. Let $\Phi(\X,\Y)=\X+\Y$ be the additive formal group law. Then $\End(\Phi)=C\{\{\tau\}\}$, the subring, under addition and composition, of $\X C[[\X]]$ of power series in the monomials $\X^{q^i}$, $i\ge0$, with coefficients in $C$. It can be shown that $\phi$ extends uniquely to a continuous ring homomorphism $\wh\phi\colon A_\Fp\to\End(\Phi), a\mapsto \wh\phi_a$. This uses that elements in $\Fp$ map to topologically nilpotent elements in $C$ under $\gamma$. This defines the structure of a formal $A_\Fp$-module $\wh\phi_\Fp=(\Phi,\wh\phi)$ on $\Phi$. Moreover the height of the formal $A_\Fp$-module $\wh\phi_\Fp\pmod{\Fm_C}$ agrees with the height of the Drinfeld $A$-module $\phi\pmod{\Fm_C}$.
\end{Ex}

Let $\breve k$ be an $\CO$-algebra via reduction, i.e., via the canonical maps $\CO\to\breve\CO\to\breve k$. Let $\overline X$ be a formal $\CO$-module over $\breve k$ of finite height $h>0$. A deformation of $\overline X$ to $C\in \CNL$ is a formal $\CO$-module $X_C$ over $C$ whose reduction modulo $\Fm_C$ is equal to $\overline X$. Two deformations $X_C$ and $X'_C$ to $C$ are isomorphic if there exists an isomorphism of formal $\CO$-modules over $C$ that reduces to the identity modulo $\Fm_C$. Since $h$ is finite, by \cite[Prop.~4.1]{DrinfeldEM} there is at most one such isomorphism.
\begin{Thm}[{\cite[Prop.~4.2]{DrinfeldEM}}]\label{Thm:UnivFGp}
The functor $\CNL\to\Sets$ that associates to $C\in\CNL$ the set of deformations of $\overline X$ to $C$ up to isomorphism is representable. The universal ring $ R_{\overline X}$ is a power series ring over $\breve\CO$ in $h-1$ indeterminates. 
\end{Thm}
\begin{Def}
The universal formal group over~$ R_{\overline X}$ is denoted by $ X_{\overline X}$.
\end{Def}

To recall the notion of $\CO$-divisible module (again of dimension $1$), we need some preparations. We fix a ring $B$ in $\CNL$. Following \cite{Taguchi}, for $R$ any ring, we define an $R$-module scheme over $B$ to be a pair $(\CG,\phi)$, where $\CG$ is a commutative group scheme over $B$ and $\phi\colon R\to \End(\CG)$ is a ring homomorphism. A map $(\CG,\phi)\to(\CG',\phi')$ of $R$-module schemes is a map $\CG\to\CG'$ of group schemes over $B$ that is equivariant for the $R$-action. 

If $\CG$ is finite flat over $B$, then one can define the \'etale and connected parts $\CG^\et$, $\CG^\loc$ of $\CG$, and one has a short exact sequence $0 \to \CG^\loc\to \CG\to \CG^\et \to 0$; see \cite[1.4]{Tate-pDivGps}. Because any endomorphism of $\CG$ preserves $\CG^\loc$, if $\CG$ carries an $R$-action, then the short exact sequence is one of $R$-module schemes. For the following, we assume that $K$ has positive characteristic. Then the field $k$ is canonically a subring of~$\CO$. 

\begin{Def}[{\cite[\S~4]{DrinfeldEM}, \cite[\S~1]{Taguchi}}]
Let $r$ be in $\BN$. An \emph{$\CO$-divisible module of rank $r$ over $B$} is an inductive system $\CF=(\CF_n,i_n)_{n\in\BN}$ such that for all $n\in\BN$ the following hold:
\begin{enumerate}
\item $\CF_n$ is a finite flat group scheme over $B$ that carries an $\CO$-module structure.
\item There is a closed immersion $\CF_n\into\BG_{a,B}$ of $k$-module schemes.\footnote{We restrict to $\CO$-modules of dimension $1$ and therefore suppress the dimension in the definition.}
\item 
The order of $\CF_n$ over $B$ is $q_K^{rn}$,
\item 
The following sequence of $\CO$-module schemes over $B$ is exact
\[
0\to \CF_n\to[i_n] \CF_{n+1}\to[\varpi_K^n]\CF_{n+1}
\]
\end{enumerate}
A \emph{morphism of $\CO$-divisible modules over $B$} is a morphism of inductive systems of $\CO$-module schemes. 
\end{Def}
Given an $\CO$-divisible module $\CF=(\CF_n)_{n\ge1}$, the connected and \'etale parts $\CF^\loc=(\CF_n^\loc)_{n\ge1}$ and $\CF^\et=(\CF_n^\et)_{n\ge1}$ form $\CO$-divisible modules as well and one has a degree-wise short exact sequence of $\CO$-divisible modules $0\to \CF^\loc\to \CF\to\CF^\et\to0$. If $\CF^\et=0$, we call $\CF$ local.

Concerning $\CF^\et$ note that since $\breve k$ is algebraically closed, one has an isomorphism of $\CO$-module schemes between $\CF_n^\et\mod{\Fm_B}$ and the constant $\CO$-module scheme $\underline{\CO^s/\varpi_K^n \CO^s }$ for $s$ the rank of $\CF^\et$. Hensel lifting shows that the same isomorphism holds over $B$. Drinfeld writes $\CF^\et= \underline{K^s/\CO^s}$.

To analyze $\CF^\loc$, we present in the following paragraphs, up to and including \autoref{Prop:LocalODiv-Equiv}, some results that are implicitly stated in \cite[\S~4]{DrinfeldEM} and are straightforward to deduce from  \cite[\S~1]{Taguchi}. Suppose that $\CF^\loc$ is non-trivial. Then $\dirlim \CF_n^\loc\cong \Spf B[[\X]]$, and this isomorphism is one of formal $k$-module schemes, if we identify $\Spf B[[\X]]$ with the formal completion of $\BG_{a,B}$ at the zero section. The action of $\CO$ on $\CF^\loc$ induces an $\CO$-action on the formal additive group over $B$. The resulting formal $\CO$-module will be denoted by~$X_\CF$.

Conversely, let $X=(\Phi,[\cdot]_X)$ be a formal $\CO$-module over $B$ whose reduction to $\breve k$ has finite height $h$. Then a local divisible $\CO$-module is defined as follows: For $n\in\BN$, write $[\varpi_K^n]_{X}=H_n u_n$ uniquely with $H_n\in B[\X]$ monic of degree $q_K^{hn}$ and $H_n\pmod{\Fm_B}=x^{q_K^{hn}}$, and $u_n\in B[[\X]]$ a unit. Then $X[\varpi_K^n]:=\Spec B[\X]/(H_n)=\Spec B[[\X]]/([\varpi_K^n]_{X})$ is a finite flat scheme over $B$ and one can verify for all $n\ge1$ that
\begin{enumerate}
\item The formal $\CO$-module structure $[\cdot]_X$ defines an $\CO$-action on $X[\varpi_K^n]$, and in such a way that the closed immersion $X[\varpi_K^n]\into \BG_{a,B}$ is one of $k$-module schemes.
\item One has a short exact sequence $0\!\to\! X[\varpi_K^{n}]\!\to\! X[\varpi_K^{n+1}]\!\stackrel{\varpi_K}\to \!X[\varpi_K^{n+1}]$ of $\CO$-module schemes.
\end{enumerate}
The resulting $\CO$-divisible local group is denoted by $\CF_X$. 
\begin{Prop}\label{Prop:LocalODiv-Equiv}
The constructions $X\mapsto \CF_X$ and $\CF\mapsto X_\CF$ define mutual inverses between the set of local divisible $\CO$-modules $\CF$ of rank $h$ and the set of formal $\CO$-modules $X=(\Phi,[\cdot]_X)$ such that $x\pmod{\Fm_B}$ has height~$h$. 
\end{Prop}

\begin{Ex}[{\cite[before Prop.~5.4]{DrinfeldEM}}]\label{Ex:DM->ODivlGp}
Let $C$, $\phi$, $\CO$, $\wh\phi$ be as in \autoref{Ex:DM->FormalGp}, and let $g_\Fp$ be a generator of the ideal $\Fp^{n_\Fp}\subset A$. Let $r$ be the rank of $\phi$ and $h$ its height. One verifies the following:
\begin{enumerate}
\item For $n\ge0$, the scheme $\phi[\Fp^{nn_\Fp}]:=C[\X]/(\phi^n_{g_\Fp}(\X))$ is finite flat over $C$ and possess an $A_\Fp$-module structure via $\phi$.
\item The sequence $(\phi[\Fp^n])_n$ with $\phi[\Fp^n]\into\phi[\Fp^{n+1}]$ given by inclusion defines a divisible $A_\Fp$-module $\phi[\Fp^\infty]$ over~$C$ of height~$r$.
\item One has an isomorphism $\CF^\loc\cong (\wh\phi_\Fp[\Fp^n])_{n\ge1}$. 
\item The rank $h$ of $\CF^\loc$ is the height of $\phi\pmod{\Fm_C}$ over $\breve k$, and one has $\CF^\et\cong(F_\Fp/A_\Fp)^{r-h}$.
\end{enumerate}
\end{Ex}

Let $\CF_{\breve k}$ be an $\CO$-divisible module of rank $r$ over $\breve k$. There is an obvious notion of a deformation $\CF_{\breve k}$ to $\CO$-divisible modules over rings $C$  in $\CNL$ and this defines a functor $\CNL\to\Sets$.
\begin{Thm}[{\cite[Prop.~4.5]{DrinfeldEM}}]\label{Thm:DrinfODivIsom}
Suppose that $\CF_{\breve k}^\loc$ has height $h>0$. Then the functor $\CNL\to\Sets$ that associates to $C\in\CNL$ the set of deformations of $\CF_{\breve k}$ to $C$ up to isomorphism is representable by some ring $ R_{\CF_{\breve k}}$ in $\CNL$. The universal ring $ R_{\CF_{\breve k}}$ is a power series ring over $\breve\CO$ in $r-1$ indeterminates.
\end{Thm}

From here on, we let $\CO:=A_\Fp$ with $\Fp$ a closed point of $\Spec A$ as in the introduction. We let $\phi_0\colon A\to \breve k\{\tau\}$ be a Drinfeld-module of rank $r$ whose characteristic is given by $A\to A_\Fp=\CO\to\breve\CO\to\breve k$, for our chosen $\breve\CO$. Let $I\subset A$ be a proper non-zero ideal with $I+\Fp=A$. Choosing a level $I$-structure for $\phi_0$, which can be done over $\breve k$, defines a point of $\FM_{I,\Fp}^d(\breve k)$ which we denote by $x$. Then $x$ defines a Drinfeld $A$-module $\phi_x$ that is
isomorphic to $\phi_0$ (over $\breve k$) together with a level $I$-structure. A deformation of $\phi_0$ is a Drinfeld $A$-module $\phi\colon A\to C\{\tau\}$, in standard form, up to isomorphism. By Hensel's Lemma, the level $I$-structure on $\phi_0$ extends uniquely to a level $I$-structure of $\phi$ over $C$. 
Hence one can identify deformations of $\phi_0$ with morphisms $\Spec C\to \FM^r_I$ that when composed with $\Spec\breve k\to\Spec C$ yield $x$. The following is Drinfeld's analog of the Serre-Tate theorem for Drinfeld $A$-modules.
\begin{Thm}[{\cite[5.C, in part.~Prop.~5.4]{DrinfeldEM}}]
The following holds
\begin{enumerate}
\item The functor $\CNL\to\Sets$ of deformations of $\phi_0$ is representable by the completion of the stalk of $\CO_{\FM_I^r\otimes_{A[1/I]}\breve\CO}$ at $x$; in particular, this completion is independent of the choice of~$I$.
\item The natural transformation from deformations of $\phi_0$ to deformations of the $\CO$-divisible group $\phi_0[\Fp^\infty]$ defined in \autoref{Ex:DM->ODivlGp}, is an isomorphism. I.e., there is a natural isomorphism of $\breve\CO$-algebras
\[  R_{\phi_0[\Fp^\infty]} \longto \wh\CO_{\FM_I^r\otimes_{A[1/I]}\breve\CO,x}.\]
\end{enumerate}
\end{Thm}

\tempnewpage
\section{The result of Strauch}
Let $K$, $\CO$, $k$, $\breve K$, $\breve\CO$, $\breve k$ and $\CNL$ be as in the previous section. Let $\overline X$ be a formal group over $k$ of height $h$ and let $R_{\overline X}$ and $X_{\overline X}$ be as in \autoref{Thm:UnivFGp}. The following is from \cite[\S~1,2]{Strauch}. First one may choose an identification $R_{\overline X}\cong \breve\CO[[u_1,\ldots, u_{h-1}]]$ such that the multiplication by $\varpi_K$ on $X_{\overline X}$ is given by a power series $[\varpi_K]_{X_{\overline X}}\in R_{\overline X}[[\X]]$ with the property that for all $i=0,\ldots,h$ one has
\begin{equation}\label{eq:FGpLcongruences}
[\varpi_K]_{X_{\overline X}} \equiv u_i\X^{q_K^i}\pmod{(u_0,\ldots,u_{i-1},\X^{q_K^i+1})},
\end{equation}
with the conventions $u_0=\varpi_K$ and $u_h=1$.

For $m\in\{0,\ldots,h-1\}$ put
\[R_{m}:=\breve \CO[[u_1,\ldots,u_{h-1}]]/(u_0,\ldots,u_{m}) \]
with the convention that $R_0 = R_{X_{\overline X}}$. Then the closed reduced subscheme of $\Spec R_0$ where the height of the connected component of $X_{\overline X}[\varpi_K^\infty]$ is at least $m$ is equal to $\Spec R_m$, and the open part of $\Spec R_{m}$ where the height of the connected component is equal to $m$ is
\[ U_m := \Spec R_m\setminus  V(u_m).\]
Let $\kappa_m$ be the field of fractions of $R_m$ and put $\eta_m=\Spec\kappa_m$. Let $\overline\kappa_m$ be an algebraic closure of $\kappa_m$ and put $\overline\eta_m=\Spec \overline\kappa_m$. Fix a positive integer $n$. Denote by
\[\Tate_{X_{\overline X},\eta_m}:=\invlim_n  X_{\overline X}[\varpi_K^n]_{\eta_m}(\overline\kappa_m) \]
the Tate-module of $X_{\overline X}$ at $\eta_m$. It is a free $\CO$-module of rank $h-m$. The absolute Galois group $\pi_1(\eta_m,\overline\eta_m)$ of $\kappa_m$ acts $\CO$-linearly on it. We denote the resulting representation by
\[\rho_{X_{\overline X},m}\colon  \pi_1(\eta_m,\overline\eta_m)\longto \Aut_{\CO}(\Tate_{X_{\overline X},\eta_m})\cong \GL_{h-m}(\CO). \]
It clearly factors via $\pi_1(U_m,\overline\eta_m)$. Then 
\cite[Thm.~2.1]{Strauch}, asserts:
\begin{Thm}\label{Thm:Strauch}
For any $m\in\{0,\ldots,h-1\}$ the homomorphism $\rho_{X_{\overline X},m}$ is surjective.
\end{Thm}

\tempnewpage

\section{Proof of \autoref{Thm:Main}}\label{sec:proof}

Let the notation be as in \autoref{Sect:One}. Set in addition $\CO=A_\Fp$, $K=\Frac\CO$, $k=A/\Fp$ and take $\breve K,\breve\CO,\breve k$ as in \autoref{Sec:Two}. Let $\xi_\Fc\in \FM_{I,\Fp}^r(\breve k)$ be a supersingular point in the component of $\FM_I^r$ labelled by $\Fc$.\footnote{Supersingular points exist; and via the action of Hecke correspondences, which preserves the supersingular locus, they can be seen to lie in every component.} Consider the following canonical morphisms of schemes
\[ \xymatrix{\Spec \wh \CO_{\FM_{I,\overline\Fp}^r,\xi_\Fc}\ar[r]\ar@/^1pc/[rr]^{\wh\iota} &\Spec \CO_{\FM_{I,\overline\Fp}^r,\xi_\Fc}\ar[r]^{\iota}&\FM_{I,\overline\Fp}^r \rlap{,}}\]
with $\wh \CO_{\FM_{I,\overline\Fp}^r,\xi_\Fc}$ the completion of the local ring $\CO_{\FM_{I,\overline\Fp}^r,\xi_\Fc}$. Denote by $\wh\eta_\Fc$ the generic point of $ \Spec \wh \CO_{\FM_{I,\overline\Fp}^r,\xi_\Fc}$ and choose a minimal geometric point $\overline{\wh\eta}_\Fc$ over $\wh\eta_\Fc$ together with a map $\overline{\wh\eta}_\Fc\to \overline\eta_\Fc$. Let $\FM_{I,\overline\Fp}^{r,\ord}\subset \FM_{I,\overline\Fp}^r$ be the locus of ordinary Drinfeld $A$-modules. It is an open subscheme since its complement is defined by the vanishing of the coefficient of $(\phi^r_{I,\overline\Fp})_{g_\Fp}\in M_{I,\overline\Fp}^r[\tau]$ in degree $\deg g_{\Fp}$, where $M_{I,\overline\Fp}^r$ is the coordinate ring of the affine scheme $\FM_{I,\overline\Fp}^r$, and hence its complement is closed in $\FM_{I,\overline\Fp}^r$. We obtain a corresponding diagram of fundamental groups with continuous group homomorphisms
\begin{equation}
\xymatrix{ 
\pi_1\big(\wh\iota^{-1}( \FM_{I,\overline\Fp}^{r,\ord} ) ,\overline{\wh\eta}_\Fc\big)\ar[r]& 
 \pi_1\big(\iota^{-1}( \FM_{I,\overline\Fp}^{r,\ord} ),\overline\eta_\Fc \big)\ar[r]&
  \pi_1(\FM_{I,\overline\Fp}^{r,\ord},\overline\eta_\Fc)\\ 
\pi_1(\wh\eta_\Fc,\overline{\wh\eta}_\Fc)\ar[u]\ar[r]&  \pi_1(\eta_\Fc,\overline\eta_\Fc)\ar[u]\ar[ur]\rlap{.}\\
}
\end{equation}
Over $\FM_{I,\overline\Fp}^{r,\ord}$ the Tate-module 
\[\Tate_\Fp\phi_{\eta_\Fc}=\invlim_n \phi_{\eta_\Fc}^r[\Fp^n](\overline \kappa_{\eta_\Fc}), \]
 is free over $\CO$ of rank $r-1$, and we have continuous homomorphisms.
\begin{equation}\label{eq:TateFromPhiAtX}
\pi_1(\wh\eta_\Fc,\overline{\wh\eta}_\Fc)\longto \pi_1(\eta_\Fc,\overline\eta_\Fc)\longto  
\pi_1(\FM_{I,\overline\Fp}^{r,\ord},\overline\eta_\Fc)  \longto \Aut_\CO(\Tate_\Fp\phi_{\eta_\Fc}).\end{equation}
By \autoref{Thm:DrinfODivIsom}, the ring $\Spec \wh \CO_{\FM_{I,\overline\Fp}^r,\xi_\Fc}$ is naturally identified with the special fiber of the universal deformation ring of the $\CO$-divisible module $\phi_{\xi_\Fc}[\Fp^\infty]$. Because $\xi_\Fc$ is supersingular, it is a local $\CO$-divisible module of rank $r$, and thus by \autoref{Prop:LocalODiv-Equiv} it arises from a formal $\CO$-module of height $r$. Now by  \autoref{Thm:Strauch} of Strauch, the composition of the maps in \eqref{eq:TateFromPhiAtX} is surjective. 

We have thus proved the following result.

\begin{Thm}\label{Thm:Main2}
The monodromy representation $\pi_1(\wh\eta_\Fc,\overline{\wh\eta}_\Fc)\longto \Aut_\CO(\Tate_\Fp\phi_{\eta_\Fc})$ is surjective.
\end{Thm}

Hence the map $\pi_1(\eta_\Fc,\overline\eta_\Fc)\to \Aut_\CO(\Tate_\Fp\phi_{\eta_\Fc})$ is surjective, as well. This completes the proof of \autoref{Thm:Main}.

\begin{Rem}\label{Rem:DecompGp}
One can think of the image of $\pi_1(\wh\eta_\Fc,\overline{\wh\eta}_\Fc)$ in $ \pi_1(\FM_{I,\overline\Fp}^{r,\ord},\overline\eta_\Fc)$ as the decomposition group at $\xi_\Fc$. From this viewpoint, \autoref{Thm:Strauch} says that already the image of this decomposition group surjects onto $\Aut_\CO(\Tate_\Fp\phi_{\eta_\Fc})\cong\GL_{r-1}(\CO)$. According to the same theorem, the decomposition groups at points of height $m<r$ map onto a natural subgroup of $\Aut_\CO(\Tate_\Fp\phi_{\eta_\Fc})$ isomorphic to $\GL_{m-1}(\CO)$.
\end{Rem}

\tempnewpage

\section{An example}\label{Sec:Example}

For the remainder of this article we specialize to the case
$A=\BF_q[t]$ and level $I=tA$. We also set $B=A[\frac{1}{t}] = \BF_q[t,\frac{1}{t}]$, we let $\Fp\in A$ be a non-zero prime (monic irreducible polynomial) and suppose $\Fp\neq t$ (otherwise, just replace $t$ by $t+1$), and we write $|\Fp| = q^{\deg(\Fp)}$. Let $\kappa_{\Fp} = A/\Fp$ with algebraic closure $\bar{\kappa}_{\Fp}$. 
As a preparation for \autoref{Application}, in the present section we will work out an explicit example of the main result.

We start by recalling Pink's explicit description of $\FM_{t}^r$ \cite{PinkSchieder,PinkCompactification}, see also \cite[Theorem 2]{BreuerExplicit} for more details:
Let $V$ be an $\BF_q$-vector space of dimension $r\geq 1$ and write $V'=V\smallsetminus\{0\}$.
Denote by $S_V=\Sym_B(V)$ the symmetric algebra of $V$ over $B$ and by $K_V$ the fraction field of~$S_V$. Denote by $RS_{V,0} = B[\frac{v}{v'} \; |\; v,v\in V']$ the subalgebra of $K_V$ generated over $B$ by quotients of non-zero elements of $V$. Then the base-change of $\FM_t^r$ to $\Spec B$ is given by
\[
\FM_{t,B}^r = \Spec RS_{V,0},
\]
which has geometrically irreducible fibres.
Furthermore, for any fixed $v_1\in V'$, the universal Drinfeld module $\phi=\phi^r_{\eta}$ on $\FM_{t,B}^r$ is determined by
\[
t \mapsto \phi_t(X) = tX\prod_{v\in V'}\left(1 - \frac{v_1}{v}X\right) \in RS_{V,0}[X]
\]
with level structure
\[
V \stackrel{\sim}{\longrightarrow} \phi[t]; \quad v \mapsto \frac{v}{v_1}.
\]
The base change of the moduli scheme $\FM_{t,B}^r$  to $\overline\kappa_\Fp$ is 
$\FM_{t,\overline{\Fp}}^r = \Spec \big( RS_{V,0}\otimes_B\overline{\kappa}_{\Fp} [X] \big)$, with universal Drinfeld module the reduction of $\phi$ modulo $\Fp$, which we denote $\overline\phi$. We have
\[
\overline\phi_{\Fp^n}(X) = \overline\phi^{\et}_{\Fp^n}(X^{|\Fp|^n}),
\]
where $\overline\phi^{\et}_{\Fp^n}(X) \in RS_{V,0}\otimes_B\overline\kappa_\Fp[X]$ is a separable $\BF_q$-linear polynomial of degree $|\Fp|^{n(r-1)}$.  
The outer terms in the local-\'etale decomposition
\[
0 \to \overline\phi[\Fp^n]^{\loc} \to \overline\phi[\Fp^n] \to \overline\phi[\Fp^n]^{\et} \to 0
\] 
are given by
\[
\overline\phi[\Fp^n]^{\loc} = \Spec \big( RS_{V,0}\otimes_B \overline{\kappa}_{\Fp}[X]/\langle X^{|\Fp|^n} \rangle \big) \qquad\text{and}\qquad
\overline\phi[\Fp^n]^{\et} = \Spec \big( RS_{V,0}\otimes_B \overline{\kappa}_{\Fp}[X]/\langle \overline\phi^\et_{\Fp^n}(X) \rangle \big).
\]

 Denote by $\overline\kappa_{\eta}$ the fraction field of $RS_{V,0}\otimes_B \overline{\kappa}_\Fp$, which is the function field of $\FM_{t,\overline{\Fp}}^r$ over $\overline\kappa_{\Fp}$, and by $\overline\kappa_{\eta}(\overline{\phi}[\Fp^n]^{\et})$ the splitting field of $\overline{\phi}^{\et}_{\Fp^n}(X)$ over $\overline\kappa_\eta$.

Now Theorem \ref{Thm:Main} says the following: For every positive integer $n$, we have
\begin{equation}\label{eq:MainPink}
\Gal\big(\overline\kappa_{\eta}(\overline{\phi}[\Fp^n]^{\et}) / \overline\kappa_\eta\big) \cong \GL_{r-1}(A/\Fp^n).
\end{equation}

\section{An application}
\label{Application}

In this last section, we consider a variant of the above example and answer a question posed in \cite{BR16}.

Suppose $g_1,g_2,\ldots,g_{r-1}$ are algebraically independent over $\BF_q(t)$ and set $L=\BF_q(t,g_1,\ldots,g_{r-1})$, a rational function field of transcendence degree $r$ over~$\BF_q$. 

We define the Drinfeld module $\psi : A \to L\{\tau\}$ by
\[
t\mapsto \psi_t(X) = tX + g_1X^q + \cdots + g_{r-1}X^{q^{r-1}} + X^{q^r} \in L[X].
\]
It is shown in \cite[Thm. 6]{BreuerExplicit} that, for every non-zero proper ideal $\Fn\subset A$, 
\[
\Gal(L(\psi[\Fn])/L) \cong \GL_r(A/\Fn).
\]
%
%
Our goal is to prove a similar result in special characteristic.

Denote by $L_t = L(\psi[t])$ the splitting field of $\psi_t(X)$ over $L$, and set
\[
RS_{t} = B\big[v,\frac{1}{v} \;|\; 0\neq v\in \psi[t]\big] \subset L_{t},
\]
the subalgebra of $L_t$ generated over $B$ by $v$ and $\frac{1}{v}$ for $0\neq v\in\psi[t]$. It is a graded ring if we set $\deg(v)=1$ for all $0 \neq v \in\psi[t]$.
We have $\psi_t(X)\in RS_t[X]$.

Fix a non-zero $t$-torsion point $0\neq v_1\in\psi[t]$, and consider the isomorphic Drinfeld module 
$\phi = v_1^{-1}\psi v_1$ over $L_t$. 
We denote by $L_{t,0}=L(\phi[t])\subset L_t$ the splitting field of $\phi_t(X)$ over $L$, and set 
\[
RS_{t,0} = B\big[\frac{v}{v'} \;|\; v,v'\in \psi[t],\; v'\neq 0\big] \subset L_{t,0}.
\]
This is the degree zero component of $RS_t$.

We have $\phi_t(X) = tX\prod_{0\neq v\in\psi[t]}\left(1-\frac{v_1}{v}X\right) \in RS_{t,0}[X]$.

By \cite[Thm. 5]{BreuerExplicit} and its proof, there is an isomorphism 
\begin{equation}\label{eq:theta}
\theta : \Spec (RS_{t,0}) \stackrel{\sim}{\longto}\FM_{t,B}^r,
\end{equation}
and $\phi$ is the pullback via $\theta$ of the universal Drinfeld module described in Section \ref{Sec:Example}. 

Now consider the reduced Drinfeld modules $\overline{\psi}$ and $\overline{\phi}$ over $RS_{t}\otimes_B\kappa_\Fp$ and $RS_{t,0}\otimes_B\kappa_\Fp$, respectively.

Again, for every positive integer $n$, we have $\overline\psi_{\Fp^n}(X) = \overline\psi^{\et}_{\Fp^n}(X^{|\Fp|^n})$, where $\overline\psi^{\et}_{\Fp^n}(X) \in RS_{t}\otimes_B \overline\kappa_{\Fp}[X]$ is a separable $\BF_q$-linear polynomial of degree $|\Fp|^{n(r-1)}$, and 
\[
\overline\psi[\Fp^n]^{\et} = \Spec \big( RS_{t}\otimes_B \overline{\kappa}_{\Fp}[X]/\langle \overline\psi^\et_{\Fp^n}(X) \rangle \big).
\]
Analoguous statements hold for $\overline\phi$. 

We define the $A$-field $\ell = \overline{\kappa}_{\Fp}(g_1,g_2,\ldots,g_{r-1})$, equipped with the homomorphism $A \to A/\Fp \subset \overline\kappa_\Fp \subset \ell$. The Drinfeld module $\overline\psi$ is defined over $\ell$. Denote by $\ell(\overline\psi[\Fp^n]^\et)$ the splitting field of $\overline\psi^\et_{\Fp^n}(X)$ over $\ell$. We have
\begin{Thm}\label{Thm:Main3}
$\Gal\big(\ell(\overline{\psi}[\Fp^n]^{\et})/\ell\big) \cong \GL_{r-1}(A/\Fp^n)$.
\end{Thm}



\begin{proof}
We define the following fields.
\begin{eqnarray*}
\ell_{t,0} & = & \overline\kappa_{\Fp}(\overline\phi[t]) = \Frac\big( RS_{t,0}\otimes_B \overline\kappa_{\Fp} \big),\\
\ell_t & = & \ell(\overline\psi[t]) = \ell_{t,0}(v_1) = \Frac\big( RS_t\otimes_B \overline\kappa_{\Fp} \big),
\end{eqnarray*}
and consider the following field extensions.
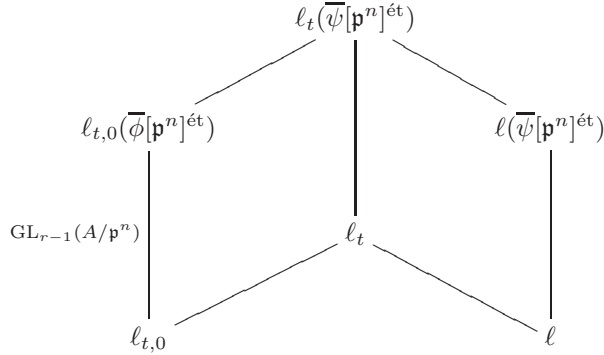
\[
\xymatrix{
& \ell_t(\overline\psi[\Fp^n]^\et) & \\
\ell_{t,0}(\overline\phi[\Fp^n]^\et)\ar@{-}[ur] & & \ell(\overline\psi[\Fp^n]^\et)\ar@{-}[ul] \\
& \ell_t\ar@{-}[uu] & \\
\ell_{t,0}\ar@{-}[uu]^{\GL_{r-1}(A/\Fp^n)}\ar@{-}[ur] & & \ell\ar@{-}[uu]\ar@{-}[ul]
}
\]

The field $\ell_{t,0}$ is (isomorphic via $\theta$ to) the function field of $\FM_{t,\overline\Fp}^r = \FM_{t,B}^r\times_{\Spec B}\Spec\overline\kappa_{\Fp}$ over $\overline\kappa_\Fp$.
It follows from (\ref{eq:MainPink}) that $\Gal\big( \ell_{t,0}(\overline\phi[\Fp^n]^\et)/\ell_{t,0} \big) \cong \GL_{r-1}(A/\Fp^n)$, and our goal is to show that the other two vertical extensions have this same Galois group.

We write
\[
\overline\phi_t(X) = \overline{t}X + c_1X^q + \cdots + c_rX^{q^r} \in RS_{t,0}\otimes_B \overline{\kappa}_\Fp[X],
\]
where $\bar{t}$ denotes the image of $t$ in $A \to A/\Fp \subset \overline\kappa_{\Fp}$, $c_i = v_1^{q^i-1}g_i$ for $i=1,\ldots,r-1$ and $c_r = v_1^{q^r-1}$. 
Because $1$ is a $t$-torsion point, we have the algebraic relation $0=\overline\phi_t(1)=\overline{t}+c_1+c_2+\cdots+c_r$.
Observe also that $\overline\kappa_\Fp\subset\ell_{t,0}$ contains roots of unity of all orders prime to $p$.
It follows that $\ell_t = \ell_{t,0}(v_1) = \ell_{t,0}(\sqrt[q^r-1]{c_r})$ is a Kummer extension of $\ell_{t,0}$.

Let $\xi\in \FM_{t,\overline\Fp}^r(\overline\kappa_{\Fp})$ correspond to a supersingular Drinfeld module 
\[
\phi^\xi_t(X) = \overline{t}X + s_1X^q + \cdots + s_rX^{q^r} \in  \overline{\kappa}_\Fp[X].
\] 
Then the completion $\widehat{\ell_{t,0}}$ of $\ell_{t,0}$ at $\xi$ contains the ring of formal power series 
\[
\overline\kappa_{\Fp}[[c_2-s_2,\cdots,c_{r}-s_{r}]].
\]  
This ring, in turn, contains 
\[
v_1 = \sqrt[q^r-1]{c_r} = \sqrt[q^r-1]{(c_r - s_r)+s_r}=\sqrt[q^r-1]{s_r}\sum_{i=0}^\infty \binom{\frac1{q^r-1}}{i} \left(\frac{c_r-s_r}{s_r} \right)^i
\] 
since $s_r\in\overline\kappa_{\Fp}^\times$ and $q^r-1$ is not divisible by the characteristic $p$.

This implies that $\overline\psi$ and $\overline\phi$ are isomorphic over $\widehat{\ell_{t,0}}$. Also, 
$\widehat{\ell_{t,0}}$ contains $\ell_t$ and 
$\widehat{\ell_{t,0}}(\overline\phi[\Fp^n]^\et) = \widehat{\ell_{t,0}}(\overline\psi[\Fp^n]^\et)$.

Theorem \ref{Thm:Main2} implies that the Galois representation 
\[
\Gal\big(\widehat{\ell_{t,0}}^\sep/\widehat{\ell_{t,0}}\big) \to \Tate_\Fp \overline\phi
\]
is surjective. Since $\overline\phi$ and $\overline\psi$ are isomorphic over $\widehat{\ell_{t,0}}$, the same holds for $\Tate_\Fp \overline\psi$, and in particular the Galois representation
\[
\Gal\big(\ell_{t}^\sep/\ell_t\big) \to \Tate_\Fp \overline\psi
\]
is surjective. This implies that 
\[
\Gal\big(\ell_t(\overline\psi[\Fp^n]^\et)/\ell_t\big) \cong \GL_{r-1}(A/\Fp^n).
\]

Finally,
we have
\[
[\ell(\overline\psi[\Fp^n]^\et) : \ell] \geq [\ell_t(\overline\psi[\Fp^n]^\et) : \ell_t] = \#\GL_{r-1}(A/\Fp^n).
\]
Since $\Gal\big(\ell(\overline\psi[\Fp^n]^\et)/\ell\big)$ is isomorphic to a subgroup of $\GL_{r-1}(A/\Fp^n)$, it must be isomorphic to the whole group. 
This completes the proof of Theorem \ref{Thm:Main3}.
\end{proof}

Finally, we address \cite[Question 4.5]{BR16}. For this, we must first recall the construction of Drinfeld modular polynomials from \cite{BR16}.

Denote by $C$ the subring of $A[g_1,\ldots,g_{r-1}]\subset L$ generated by monomials of the form $ag_1^{e_1}\cdots g_{r-1}^{e_{r-1}}$ satisfying $a\in A$ and $\sum_{k=1}^{r-1}e_k(q^k-1) \equiv 0 \bmod q^r-1$. Then the elements of $C$ are the isomorphism invariants of rank $r$ Drinfeld $A$-modules, i.e. $\Spec C$ is the coarse moduli scheme of Drinfeld modules of rank $r$ and no level structure, see \cite[Prop. 1.1]{BR16}.

Let $1\leq s \leq r-1$. An isogeny $f : \psi \to \psi^{(f)}$ is said to have type $(A/\Fp)^s$ if $\ker f(\overline{L}) \cong (A/\Fp)^s$, and such an isogeny is called {\em special} if $\ker f$ contains 
$U_0:=\ker \big(\psi[\Fp](\overline{L}) \to \overline\psi[\Fp](\overline{\ell})\big)$.
Because $\psi$ has ordinary reduction~$\bar\psi$, and so $U_0\cong A/\Fp$, $f$ is special if and only if its reduction is inseparable.

To each invariant $J\in C$ we associate the Drinfeld modular polynomial of type $(A/\Fp)^s$, defined by
\[
\Phi_{J,(A/\Fp)^s}(X) = \prod_{\text{$f : \psi \to \psi^{(f)}$ of type $(A/\Fp)^s$}}\big(X - J(\psi^{(f)})\big) \in C[X].
\]
This is irreducible over $L$ if its roots in $\overline{L}$ are distinct (there always exist $J\in C$ for which the roots are distinct). 

Modulo $\Fp$, we have the Kronecker congruence relation \cite[Thm. 4.4]{BR16}:
\begin{equation}\label{eq:Kronecker}
\Phi_{J,(A/\Fp)^s}(X) \equiv \Phi^{\spec}_{J,(A/\Fp)^s}(X)\cdot \big( \Phi^{\spec}_{J,(A/\Fp)^{s+1}}(X^{|\Fp|})\big)^{|\Fp|^{s-1}} \bmod \Fp,
\end{equation}
where 
\[
\Phi^{\spec}_{J,(A/\Fp)^s}(X) := \prod_{\text{$f : \psi \to \psi^{(f)}$ special of type $(A/\Fp)^s$}}\big(X - (J(\psi^{(f)}) \bmod \Fp) \big) \in \ell[X].
\]

We answer \cite[Question 4.5]{BR16} in the affirmative, as follows.

\begin{Prop}\label{Prop:OnBR16}
Suppose $J\in C$ is such that the roots of $\Phi^{\spec}_{J,(A/\Fp)^s}(X)$ in $\overline\ell$ are distinct.
Then $\Phi^{\spec}_{J,(A/\Fp)^s}(X)\in\ell[X]$ is irreducible.
\end{Prop}

\begin{proof}
If $s=1$, then $\Phi^{\spec}_{J,(A/\Fp)}(X)=X-J^{|\Fp|}$ by \cite[Example 5.1]{BR16}, and we are done. 

Now suppose that $s > 1$.
Let $R$ be the integral closure of $A[g_1,g_2,\ldots,g_{r-1}]$ in $L(\psi[\Fp])$, then $\psi[\Fp]\subset R$.
Let $f : \phi \to \phi^{(f)}$ be a special isogeny of type $(A/\Fp)^s$. Then $f(X) \in R[X]$ is an $\BF_q$-linear polynomial, and $f(X) \equiv f^{\et}(X^{|\Fp|}) \bmod \Fp$, where $f^{\et}(X) \in R\otimes_A \kappa_\Fp[X]$ is separable and $\ker f^{\et}$ is an $A$-submodule of $\overline\psi[\Fp]^{\et}$ isomorphic to $(A/\Fp)^{s-1}$.

By Theorem \ref{Thm:Main3} $\Gal(\ell^\sep/\ell)$ acts transitively on the set of such submodules of $\overline\psi[\Fp]^{\et}$, and thus also on the set of Drinfeld modules $\overline\psi^{(f)}$.
Because $J$ maps different $\overline\psi^{(f)}$ to different elements of $\overline \ell$, the group $\Gal(\ell^\sep/\ell)$, in turn, acts transitively on the roots of $\Phi^{\spec}_{J,(A/\Fp)^s}(X)$. 
\end{proof}

\begin{Rem}
Equation (\ref{eq:Kronecker}) thus describes the decomposition of the Hecke correspondence associated to $(A/\Fp)^s$-isogenies on $\FM^r_{I,\overline\Fp}$ into irreducible components with multiplicities.
\end{Rem}

\section*{Acknowledgements}
The authors are grateful to Judith Ludwig for helpful discussions and comments.
The second author would like to thank the University of Heidelberg for its hospitality and the Alexander-von-Humboldt Foundation for financial support.

\bibliographystyle{hep}
\addcontentsline{toc}{section}{References}
\bibliography{GM-DMV}

\begin{center}
\rule{8cm}{0.01cm}
\end{center}

\begin{minipage}[t]{7cm}{\small
Computational Arithmetic Geometry, IWR \\
University of Heidelberg,\\ 
Heidelberg\\ 
Germany\\
gebhard.boeckle@iwr.uni-heidelberg.de
}
\end{minipage}\hfill
\begin{minipage}[t]{7cm}{\small
School of Mathematical and Physical Sciences \\
University of Newcastle \\
Newcastle\\ 
Australia \\
florian.breuer@newcastle.edu.au}
\end{minipage}

\end{document}